\documentclass[12pt]{amsart}
\input xy
\xyoption{all}
\newcommand{\bP}{{\mathbb P}}
\newcommand{\bQ}{{\mathbb Q}}

\newcommand{\bZ}{{\mathbb Z}}

\newcommand{\cA}{{\mathcal A}}

\newcommand{\cO}{{\mathcal O}}

\newcommand{\cR}{{\mathcal R}}

\newcommand{\ox}{\overline{x}}
\newcommand{\oz}{\overline{z}}

\newcommand{\tC}{\widetilde{C}}

\newcommand{\tD}{\widetilde{D}}

\newcommand{\tX}{\widetilde{X}}

\newcommand{\ra}{\rightarrow}
\newcommand{\lra}{\longrightarrow}
\newcommand{\lla}{\longleftarrow}

\newcommand{\T}{\Theta}

\DeclareMathOperator{\Div}{{Div}}
\DeclareMathOperator{\Prym}{{Prym}}

\DeclareMathOperator{\cliff}{{cliff}}
\DeclareMathOperator{\gon}{{gon}}

\newtheorem{theorem}{Theorem}[section]
\newtheorem{lemma}[theorem]{Lemma}
\newtheorem{proposition}[theorem]{Proposition}
\newtheorem{corollary}[theorem]{Corollary}
\newtheorem{remark}[theorem]{Remark}

\numberwithin{equation}{section}

\begin{document}

\title{Counter-examples of high Clifford index to Prym-Torelli}

\author{E. Izadi}

\address{Department of Mathematics, Boyd
Graduate Studies Research Center, University of Georgia, Athens, GA
30602-7403, USA}

\email{izadi@math.uga.edu}

\author{H. Lange}

\address{Department Mathematik, Universit\"at Erlangen-N\"urnberg, Bismarckstrasse 1 1/2, D-91054 Erlangen, Germany}

\email{lange@mi.uni-erlangen.de}

\thanks{The first author was partially supported by the National Security Agency. Any opinions,
findings and conclusions or recomendations expressed in this material
are those of the author and do not necessarily reflect the views of
the National Security Agency. Part of this work was done while H. Lange was visiting the
University of Georgia: he would like to thank the University of Georgia for its hospitality.}




\maketitle

\section*{Introduction}

To any (non-trivial) \'etale double covering $\pi: \tX \ra X$ of a smooth projective curve $X$ of genus 
$g \geq 2$ one can associate a principally polarized abelian variety $P(\pi)$ of dimension $g-1$ in a canonical way,
the {\it Prym variety} of $\pi$. This induces a morphism
$$
pr_g: \cR_g(2) \ra \cA_{g-1}
$$
from the moduli space $\cR_g(2)$ of (non-trivial) \'etale double coverings of curves of genus $g$ to the moduli space
$\cA_{g-1}$ of principally polarized abelian varieties of dimension $g-1$, called the {\it Prym map}. It was shown independently
by Kanev, Friedman-Smith, Welters and Debarre, that $pr_g$ is generically injective for $g \geq 7$. On the other hand, 
Beauville remarked in \cite{beauville772} that $pr_g$ is not injective for $g \leq 10$. In \cite{Donagi81} Donagi gave 
a construction showing that $pr_g$ is not injective at any \'etale double cover of a curve $X$ admitting a map $X \ra \bP^1$
of degree 4 under some generality assumptions. Moreover, he conjectured 
(see \cite[Conjecture 4.1]{Donagi81} or \cite[p. 253]{sh}) 
that $pr_g$ is injective at any $\pi: \tX \ra X$, whenever $X$ does not admit a $g_4^1$. 

Verra showed in \cite{v} that $pr_{10}$ is not injective at any \'etale double cover of 
a general plane sectic. However, the curves which either admit a $g_4^1$ or are
plane sextics (more precisely, admit a $g^2_6$)
are exactly the curves of Clifford index $\leq 2$. So one might ask whether $pr_g$ is injective at $\pi:\tX \ra X$ whenever  
$X$ is of Clifford index $\geq 3$. It is the aim of this paper to show that this is not the case. Our main result is the
following theorem.

\begin{theorem} \label{thm01}
For any integer $N$ there is a curve $X$ of Clifford index at least $N$ such that the Prym map is not injective at any 
\'etale double cover of $X$.  
\end{theorem} 

If $\rho_4: X \ra Y$ is a ramified 4-fold cover of smooth projective curves, the tetragonal construction generalizes immediately
to associate to any \'etale double cover $\kappa: \tX \ra X$ two other \'etale double covers $\tau_i: \tC_i \ra C_i$ 
where $C_i$ is a 4-fold cover of $Y$ of the same genus as $X$. It was shown in
\cite[Paragraph 6.5]{I19} that the corresponding
Prym varieties $P = P(\kappa)$ and $P_i = P(\tau_i)$ are isogenous. In the special case $Y
\cong \bP^1$ Donagi showed that they are isomorphic. We show that, under small generality assumptions, they are also isomorphic for an arbitrary curve $Y$,
which leads to Theorem \ref{thm01}.

We know two existing proofs of Donagi's theorem, one via degeneration which is given in
\cite{Donagi92} and the other using the cohomology 
class of the Abel-Prym curve (see \cite{bl}). Neither existing proof seems to generalize to arbitrary base curves $Y$.
We give a proof of this using an explicit correspondence that induces the isogeny $P \ra P_i$. We show that this 
isogeny has degree $2^{ 2 g_X -2}$ and factors through multiplication by 2, which gives the isomorphism. This proof works in particular in the case $Y = \bP^1$,
thus giving a third proof of Donagi's theorem.

In a different direction, we consider the question of the existence of irreducible
curves representing multiples of the minimal class in Prym varieties. The
minimal cohomology class for curves in an abelian variety $A$ of dimension $g$ with
principal polarization $\T$ is
\[
\frac{[\T]^{g-1}}{(g-1)!}.
\]
Our construction above of the curves $\tC_i$ can also be done if we replace $4$-fold
covers $\rho_4$ by $n$-fold covers $\rho_n$, for any $n\geq 3$. We prove the following (see Theorem
\ref{thmclass} below).

\begin{theorem} \label{thm02}
For any unramified double cover $\kappa: \tX \ra X$ and any simply ramified $n$-sheeted cover $\rho_n: X \ra Y$ as in 
Section 1 the class of the image of $\tC_i$ in the Prym variety $P$ of $\kappa$ is 
\[
[\tC_i] = 2^{n - 1}\frac{[\T_{P}]^{g_X -2}}{(g_X -2)!}
\]
where $\T_P$ is the principal polarization of $P$ as the Prym variety of the double cover $\kappa$.
\end{theorem}
Our original motivation for wanting to produce examples of curves with highly
split jacobians (see \cite{I19}: the curves $\tC_i$ have this property) was to find examples of irreducible curves
representing multiples of the minimal class in principally polarized abelian varieties.

We work over the field of complex numbers.

\section{\bf Summary of previous results}\label{sectionsumm}

\subsection{}

Let us recall the set up and some of the results of \cite{I19}. Consider
the following maps of smooth projective curves
\[
\tX\stackrel{\kappa}{\lra} X\stackrel{\rho_n}{\lra }Y
\]
where $Y$ is of genus $g_Y$, $\rho_n$ is a {\em simply} ramified cover of
degree $n\geq 3$, $X$ has genus $g_X$ and $\tX$ is an \'etale double cover
of $X$ which is NOT obtained by base change from a double cover of
$Y$.

Then $Y$ embeds into the symmetric power $X^{ (n) }$ via the map
sending a point $y$ of $Y$ to the divisor obtained as the sum of its
preimages in $X$. Let $\tC\subset\tX^{(n)}$ be the curve defined by the
fiber product diagram

\begin{eqnarray} \label{eq1}
\tC & \lra & \tX^{ (n) } \nonumber \\
\downarrow & & \downarrow \kappa^{ (n) } \\
Y & \lra & X^{ (n) }. \nonumber 
\end{eqnarray}

In other words, the curve $\tC$ parametrizes the liftings of points of
$Y$ to $\tX$. In \cite[Lemma 1.1]{I19} we proved that $\tC$ is smooth.
As in \cite{I19} we shall assume that $\tC$
has two connected components which we denote $\tC_1$ and $\tC_2$.
This is the case for instance when the monodromy of the cover $\tX \ra Y$ factors through the Weyl group of $D_n$. In particular, this is always the case if $Y = \bP^1$.

The curve $\tC$ has an involution $\sigma$ defined as follows. Let

\[
\oz :=\ox_1 +\ldots + \ox_n
\]
be the sum of the points in a fiber of
$\rho_n$, and, for each $i$, let $x_i$ and $x_i'$ be the two preimages
of $\ox_i$ in $\tX$. Then
\[
z := x_1 +\ldots + x_n
\]
is a point of $\tC$
and
\[
\sigma (z) = x_1' +\ldots + x_n'.
\]
Let $C$ be the quotient of $\tC$ by
$\sigma$.

The degrees of the maps $\tC\ra Y$ and $C\ra Y$ are $2^n$ and $2^{ n-1
}$ respectively. It is easily seen that $\sigma$ is fixed-point-free
if $n\geq 3$. Also, we can see that for each ramification point $\ox_1
=\ox_2$ of $\rho_n$ there are $2^{n-2}$ ramification points in a fiber
of $\tC\ra Y$ obtained as $x_1 + x_1' + D_{ n-2 }$ where $D_{ n-2 }$
is one of the $2^{ n-2 }$ divisors on $\tX$ lifting $\ox_3 +\ldots
+\ox_n$.

Two liftings of $\oz$ are in the same connected component of $\tC$ if
and only if they differ by an even number of points of $\tX$. Half of
the divisors $x_1 + x_1' + D_{ n-2 }$ lie in $\tC_1$ and the other
half lies in $\tC_2$. So the degree of the map $\tC_i \ra Y$ is $2^{n-1}$ for $i=1$ and $2$,
and $\tC_1$ and $\tC_2$ have the same genus.

Writing the degree of the ramification divisor of $\rho_n$
as
\[
\deg (R_{ X/Y } ) = 2g_X -2 - n( 2g_Y -2 ),
\]
the genus of $\tC_1$ and $\tC_2$ is
\begin{equation} \label{genus(tC)}
g_{\tC_i } = 2^{ n-3 } \left( g_X -1 - (n-4 )(g_Y -1 ) \right) +1.
\end{equation}

If $n$ is odd, the involution $\sigma$ exchanges the two components of 
$\tC$, hence induces isomorphisms
\[
\tC_1\cong\tC_2\cong C.
\]

If $n$ is even, the involution $\sigma$ acts on each component of $\tC$ 
hence $C$ also has two connected components, say $C_1$ and $C_2$. For 
$n\geq 4$, since $\sigma$ is fixed-point-free, we compute the genus of 
$C_1$ and $C_2$ to be
\begin{equation} \label{genus(C)}
g_{ C_i } = 2^{ n-4 } \left( g_X -1 - (n-4 )(g_Y -1 ) \right) +1.
\end{equation}
So we have the following diagrams
\[
n\: odd\hskip160pt n\: even
\]
\[
\xymatrix{
\tC_1\cup\tC_2 \ar[d]_{\tau_1 \cup \tau_2} & \tX \ar[d]^{\kappa} \\
C =\tC_1 =\tC_2 \ar[dr]_{\mu} & X \ar[d]^{\rho_n} \\
 & Y}
\hskip70pt
\xymatrix{
\tC_1 \ar[d]_{\tau_1} & \tX \ar[d]^{\kappa} & \tC_2 \ar[d]^{\tau_2}\\
C_1 \ar[dr]_{\mu_1} & X \ar[d]^{\rho_n} & C_2 \ar[dl]^{\mu_2}\\
 & Y.}
\]

In the sequel we work with one of the components $\tC_i$, say $\tC_1$.

\section{\bf The Correspondences $S$ and $S^t$} \label{subsecXtCi}

For our constructions we shall need the following correspondences. We define $S$ to be the reduced curve 
$$
S = \{ (z=x_1 + \cdots + x_n,x) \in \tC_1 \times \tX \;|\; x = x_i \; \mbox{for some} \; i \}
$$
and $S^t$ to be its transpose, i.e.
$$
S^t = \{ (x, z=x_1 + \cdots + x_n) \in \tX \times \tC_1 \;|\; x = x_i \; \mbox{for some} \; i \}
$$
As maps of curves to divisors, $S$ and $S^t$ are given by
$$
S: \left\{\begin{array}{ccc}
          \tC_1 & \ra & \Div^n(\tX)\\
          z = x_1 + \cdots + x_n& \mapsto & x_1 + \cdots + x_n
          \end{array} \right.
$$           
and
$$
S^t: \left\{\begin{array}{ccc}
          \tX & \ra & \Div^{2^{n-2}}(\tC_1)\\
          x & \mapsto & \sum  x + x_2^{\epsilon_2} + \cdots + x_n^{\epsilon_n}
          \end{array} \right.
$$  
where $\kappa(x_i)= \overline{x}_i,\; \kappa(x)= \overline{x}, \; 
\rho_n^{-1}\rho_n(\overline{x}) = \{\overline{x},\overline{x}_2,\cdots,\overline{x}_n \}$ and
the sum is to be taken over all $\epsilon_i = 1$ or $-1$ with $\sum \epsilon_i \equiv 0
[2]$ and $x_i^1 = x_i$, $x_i^{-1} = x_i'$ for $i=2,\cdots,n$. In other words, the sum 
is taken over all divisors of $\tau_1^{-1}\mu_1^{-1}\rho_n\kappa(x)$ in $\tC_1$ containing $x$.

As in \cite{I19}, for each $k\in\{ 1,\ldots
,n\}$ and $z= x_1 + \ldots + x_n\in\tC$, we denote by
\[
[k + (n-k )'](z)
\]
the sum of all the points of $\tC$ where
$k$ of the $x_i$ are added to $(n-k)$ of the $x_i'$, the indices $i$
being all distinct. When $k$ is even, this induces a map, also denoted $[k + (n-k)']$, from $\tC_1$ to $Div (\tC_1)$.

We need the following.

\begin{lemma}\label{lembetay}

We have
\[
S^t S = \sum_{ i=0 }^{[\frac{n-1}{2}] } (n-2i) [(2i)' +
(n-2i)].
\]

\end{lemma}

\begin{proof}

Choose a point $z \in \tC_1$. Then $z=x_1 +\ldots + x_n\in\tX^{ (n) }$,
i.e., the image of $z$ in $Div^n\tX$ is $x_1 +\ldots +
x_n$. The map $S^t$ will send this to
\[
\sum_{ i=1 }^n (x_i +\tX^{ (n-1)})\cap \tC_1.
\]
It is easy to see that this is equal to the expression in the
statement of the lemma.
\end{proof}

\begin{lemma}\label{lemSSt}

For all $x\in \tX$, we have
\[
S ( S^t (x)) = 2^{ n-2 } x + 2^{n-3}\sum_{i=2}^n (x_i + x_i').
\]

\end{lemma}

\begin{proof}
Immediate from the definitions.
\end{proof}

The correspondences $S$ and $S^t$ induce homomorphisms of the corresponding jacobians in the usual way. We denote these homomorphisms by
\[
s : J\tC_1 \lra J\tX \qquad \hbox{ and } \qquad s^t : J \tX \lra J \tC_1.
\]

Assume that $n$ is even. It follows from the definition of the involution on $\tC_1$ that, on the jacobians,
\[
s\sigma = {}' s \qquad \hbox{ and } \qquad s^t {}' = \sigma s^t.
\]
Therefore the homomorphisms $s$ and $s^t$ induce homorphisms of the Prym varieties
\begin{equation} \label{e2.1}
P := \Prym(\kappa) \quad  \mbox{and}  \quad P_1 := \Prym(\tau_1)
\end{equation}
of the involutions $\sigma$ and ${}'$, which we again denote by $s$ and $s^t$:
\[
s : P_1 \lra P \qquad s^t : P \lra P_1.
\]
Note that, for $k$ even, the endomorphism induced by $[k' + (n-k)]$ on $J\tC_1$ induces an endomorphism on $P_1$. An immediate consequence of lemmas \ref{lembetay} and \ref{lemSSt} is the following.

\begin{corollary}\label{corsststs}

Suppose $n$ is even. The homomorphisms $s$ and $s^t$ between the Prym varieties $P_1$ and $P$ satisfy the identities
\[
s^t s = \sum_{ i=0 }^{[\frac{n-1}{2}] } (n-2i) [(2i)' +
(n-2i)]
\]
and
\[
s s^t = 2^{n-2}\cdot 1_{P}
\]
where $[(2i)' + (n-2i)]$ denotes the endomorphism induced on $P_1$ by $[(2i)' + (n-2i)]$ on $J\tC_1$.

In particular, $s^t$ is an isogeny from $P$ to an abelian subvariety of $P_1$.

\end{corollary}

\section{\bf Comparison of $\Prym(\tau_1)$ and $\Prym(\kappa)$ for $n=4$}

When $n=4$, from \eqref{genus(tC)} and \eqref{genus(C)} we see that
$$
\dim P_1 = \dim P  = g_X -1. 
$$

Hence, when $n=4$, Corollary \ref{corsststs} becomes

\begin{corollary}

Suppose $n=4$. The homomorphisms $s$ and $s^t$ between the Prym varieties $P_1$ and $P$ are isogenies satisfying the identities
\[
s^t s = 4\cdot 1_{P_1}
\]
and
\[
s s^t = 4\cdot 1_P.
\]

\end{corollary}

\begin{proof}

For $n=4$, the first identity in Corollary \ref{corsststs} is
\[
s^t s = 4\cdot 1_{P_1} + [2' + 2].
\]
Since $P_1$ is the image of the endomorphism $1_{J\tC_1} - \sigma$ of $J\tC_1$, and $[2'+2]\sigma = [2+2'] = [2'+2]$, we see that $[2'+2]$ induces the zero endomorphism $0_{P_1}$ on $P_1$. \\
The second identity is immediate. 
\end{proof}

\begin{corollary} \label{cor2.4} 
$$
\deg(s: P_1 \ra P) = \deg(s^t: P \ra P_1) = 2^{2 \dim P_1} = \deg(2_{P_1}) = \deg(2_P).
$$
\end{corollary}

\begin{proof}
The fact that $S^t$ is the transposed correspondence of $S$ implies that $s^t: P \ra P_1$ is the transposed endomorphism
of $s: P_1 \ra P$ with respect to the canonical principal polarizations. In particular $\deg s^t = \deg s$. So 
Corollary \ref{corsststs} implies the assertion.
\end{proof}
           
\begin{proposition} \label{prop2.5}
The isogeny $s:P_1 \ra P$ factors via the multiplication by $2$ endomorphism $2_{P_1}: P_1 \ra P_1$. 
\end{proposition} 

\begin{proof}
Recall that $J\tC_1 = H^0(\omega_{\tC_1})^*/H_1(\tC_1, \bZ)$ and $J\tX = H^0(\omega_{\tX})^*/H_1(\tX,\bZ)$.
If, by a slight abuse of notation, $s: H^0(\omega_{\tC_1})^* \ra H^0(\omega_{\tX})^*$ 
and $s^t: H^0(\omega_{\tX})^* \ra H^0(\omega_{\tC_1})^*$ 
also denote the liftings of the homomorphisms $s$ and $s^t$, 
it is a standard fact that they satisfy the relations
$$
(\alpha, s(\beta))_{\tX} = (s^t(\alpha),\beta)_{\tC_1}
$$
for all $\alpha \in H_1(\tX, \bZ)$ and $\beta \in H_1(\tC_1,\bZ)$, where $( \;,\,)$ denote the intersection products. 
Moreover the canonical principal polarizations $E_{J\tX}$ and $E_{J\tC_1}$ are given by 
$$
E_{J\tX}(v,w) = -(v,w)_{\tX}
$$
for all $v,w \in H_1(\tX,\bZ)$ and
$$
E_{J\tC_1}(v_1,w_1) = -(v_1,w_1)_{\tC_1}
$$
for all $v_1, w_1 \in H_1(\tC_1, \bZ)$. This implies
$$
E_{J\tX}(v,s(w_1)) = E_{J\tC_1}(s^t(v),w_1)
$$
for all $v \in H_1(\tX, \bZ),\; w_1 \in H_1(\tC_1,\bZ)$.

Now let 
$$
P = V^-/\Lambda^- \quad \mbox{and} \quad  P_1 = V_1^-/\Lambda_1^-
$$
where as usual $V^-$ and $V_1^-$ are the anti-invariant vector subspaces of $H^0(\omega_{\tX})^*$ 
and $H^0(\omega_{\tC_1})^*$ with respect to the liftings of the involutions $'$ and $\sigma$.
The canonical principal polarizations of $J\tX$ and $J\tC_1$ restrict to twice principal polarizations $\Xi$ on $P$
and $\Xi_1$ and $P_1$. Hence $E_{J\tX} = 2 E_{\Xi}$ and $E_{J\tC_1} = 2 E_{\Xi_1}$ and we get
$$
E_{\Xi}(v,s(w_1)) = E_{\Xi_1}(s^t(v),w_1)
$$ 
for all $v \in \Lambda^-$ and $w_1 \in \Lambda_1^-$.
Using Corollary \ref{corsststs} this implies for all $v_1, w_1 \in \Lambda_1^-$,
\begin{eqnarray*}
s^*E_{\Xi}(v_1,w_1) = E_{\Xi}(s(v_1),s(w_1)) & = & E_{\Xi_1}(s^t s(v_1), w_1)\\
& = & E_{\Xi_1}(4 v_1,w_1)\\
& = & E_{\Xi_1}(2 v_1, 2 w_1) = 2^*E_{\Xi_1}(v_1,w_1).
\end{eqnarray*}
So $s^*E_{\Xi} = 2^*E_{\Xi_1}$ and $s^* E_{\Xi}$ takes integer values on the lattice $\frac{1}{2}\lambda_1^-$. Hence the points of order $2$ map to zero by $s$ because $\Xi$ is a principal polarization and $s$ is an isogeny. In particular, by e.g. \cite[Corollary 2.4.4 page 36]{bl}, the isogeny $s$ factors via the 2-multiplication $2_{P_1}$.
\end{proof} 

As an immediate consequence we get the main result of this section.

\begin{theorem} \label{thm2.6}
Let $\rho_4:X \ra Y$ be a simply ramified cover of degree $4$ of smooth projective curves and $\kappa: \tX \ra X$ be an \'etale doble cover 
which is not obtained by base change from a double cover of $Y$. Assume that the curve $\tC$ obtained by diagram \eqref{eq1}
has $2$ connected components one of which is $\tC_1$. If $P$ and $P_1$ denote the associated Prym varieties 
defined in \eqref{e2.1}, then we have:
 
The homomorphism $s: J\tC_1 \ra J\tX$ induces an isomorphism of principally polarized abelian varieties
$$
(P_1,\Xi_1) \simeq (P,\Xi) .
$$
\end{theorem} 
\begin{proof}
According to Proposition \ref{prop2.5} the isogeny $s: P_1 \ra P$ factors as follows
\[
\xymatrix{
P_1 \ar[dr]_{2_{P_1}}  \ar[rr]^{s} & & P  \\
& P_1 \ar[ur]_{\psi} }
.\]
According to Corollary \ref{cor2.4}, $\deg (s: P_1 \ra P) = \deg 2_{P_1}$. Hence $\psi: P_1 \ra P$ is an isomorphism. By construction 
it respects the polarizations.   
\end{proof}

\section{\bf Non-isomorphy of the coverings}

In order to show that Theorem \ref{thm2.6} gives examples of the non-injectivity of the
Prym map, we have to show that for general coverings $\rho_4$ 
the coverings $\kappa$ and $\tau_1$ are 
non-isomorphic and that the Clifford index of the curve $X$ will be at least 3. First we show more generally,

\begin{lemma}

Suppose that the cover $X\stackrel{\rho_n}{\lra }Y$ has exactly two ramification points of index $1$ in one fiber and is otherwise simply ramified, then, for any choice of \'etale  double cover $\tX\stackrel{\kappa}{\lra} X$, the curve $\tC$ has $2^{n-4}$ singular points that are described as follows.

Let $y$ be the point of $Y$ such that there are two ramification points of index $1$, say $\ox_1$ and $\ox_2$ in $\rho_n^{-1} (y)$. Let $\ox_3, \ldots, \ox_{n-2}$ be the other (distinct) points in $\rho_n^{-1} (y)$. Then the singular points of $\tC$ are the points $x_1 +x_1'+ x_2 + x_2'+ x_3 + \cdots + x_{ n-2}$ where $\kappa^{-1} (\ox_1) = \{ x_1, x_1'\}$, $\kappa^{-1} (\ox_2) = \{ x_2, x_2' \}$ and, for $i\geq 3$, $x_i$ is a point of $\tX$ lying above $\ox_i$.

\end{lemma}

\begin{proof}
As in the proof of Lemma 1.1 on page 186 of \cite{I19}, we can see that $C_1$ and $C_2$
are smooth above all non-branch points of $Y$ and all simple branch points of $Y$. For all
$i$, let $x_i$ and $x_i'$ be the points of $\tX$ above $\ox_i$. We need to analyze the local
structure of $\tC$ at the points $2x_1+ 2x_2 + x_3 + \cdots + x_{ n-2}$, $x_1 + x_1'+ 2x_2
+ x_3 + \cdots + x_{ n-2}$ and $x_1 +x_1'+ x_2 + x_2'+ x_3 + \cdots + x_{ n-2}$, the cases
of the other points of $\tC$ being similar either to these or to cases we considered in
Lemma 1.1 of loc. cit.. All these points lie above the point $2\ox_1+2\ox_2+\ox_3+ \ldots
+ \ox_{n-2}$ of $X^{(n)}$ which, with our conventions, we identify with the point $y$ of
$Y$. As in that proof, since $\tC$ is defined by the fiber product diagram \eqref{eq1},
its tangent space is the pull-back of the tangent space of $Y$. The tangent space to $Y$
is a subspace of the tangent space of $X^{(n)}$ which, at the point $2\ox_1+ 2\ox_2 +
\ox_3 + \cdots + \ox_{ n-2}$, can be canonically identified with
\[
\cO_{2\ox_1} (2\ox_1) \oplus \cO_{2\ox_2} (2\ox_2) \oplus_{i=3}^{n-2}\cO_{\ox_i }(\ox_i).
\]

At the point $2x_1+ 2x_2 + x_3 + \cdots + x_{ n-2}$, the tangent space to $\tX^{(n)}$ can be canonically identified with
\[
\cO_{2x_1} (2x_1) \oplus \cO_{2x_2} (2x_2) \oplus_{i=3}^{n-2}\cO_{x_i }(x_i).
\]
The differential of $\kappa^{(n)}$ sends $\cO_{x_i} (x_i)$
isomorphically to $\cO_{\ox_i }(\ox_i)$ and sends $\cO_{2x_1} (2x_1)$ and
$\cO_{2x_2} (2x_2)$ isomorphically to $\cO_{2\ox_1} (2\ox_1)$ and $\cO_{2\ox_2} (2\ox_2)$ respectively. Hence the differential of $\kappa^{(n)}$ is an isomorphism and $\tC$ is smooth at $2x_1+ 2x_2 + x_3 + \cdots + x_{ n-2}$.

The tangent space to $\tX^{(n)}$ at $x_1+x_1' +2x_2+x_3+\ldots
x_{n-2}\in \tX^{ (n)}$ can be canonically identified with
\[
\cO_{x_1} (x_1) \oplus \cO_{x_1'} (x_1') \oplus \cO_{2x_2} (2x_2) \oplus_{i=3}^{n-2}\cO_{x_i }(x_i).
\]
The differential of $\kappa^{(n)}$ sends $\cO_{x_i} (x_i)$ and $\cO_{2x_2} (2x_2)$ isomorphically to 
$\cO_{\ox_i }(\ox_i)$ and $\cO_{2\ox_2} (2\ox_2)$ respectively, and sends $\cO_x (x)$ and
$\cO_{x' }(x')$ both isomorphically to the subspace $\cO_{\ox}(\ox)$
of $\cO_{2\ox} (2\ox)$. This case is therefore entirely similar to the case of a simple ramification point considered in the proof of Lemma 1.1 of \cite{I19} and $\tC$ is smooth at such a point.

Finally, the tangent space to $\tX^{(n)}$ at $x_1+x_1' +x_2+ x_2'+x_3+\ldots
x_{n-2}\in \tX^{ (n)}$ can be canonically identified with
\[
\cO_{x_1} (x_1) \oplus \cO_{x_1'} (x_1') \oplus \cO_{x_2} (x_2) \oplus \cO_{x_2'} (x_2') \oplus_{i=3}^{n-2}\cO_{x_i }(x_i).
\]
For all $i$, the differential of $\kappa^{(n)}$ sends $\cO_{x_i} (x_i)$ isomorphically to 
$\cO_{\ox_i }(\ox_i)$ and, for $i=1$ or $2$, sends $\cO_{x_i} (x_i)$ and
$\cO_{x_i' }(x_i')$ both isomorphically to the subspace $\cO_{\ox_i}(\ox_i)$
of $\cO_{2\ox_i} (2\ox_i)$. Its kernel is therefore two-dimensional and it follows that $\tC$ is singular at $x_1+x_1' +x_2+ x_2'+x_3+\ldots
x_{n-2}$.
\end{proof}

\begin{corollary}

If $n=4$, for a generic choice of
\[
X\stackrel{\rho_4}{\lra }Y,
\]
and any double cover
\[
\tX\stackrel{\kappa}{\lra} X,
\]
the curves $C_1$, $C_2$ and the curves $\tC_1$, $\tC_2$ are non isomorphic and non-isomorphic to $X$, respectively, $\tX$.

\end{corollary}

\begin{proof}
In the situation of the lemma, the curve $\tC$ has one singular point. Hence, one of $\tC_1$ or $\tC_2$, say $\tC_1$, is singular and the other is smooth. Since, by the description in the lemma, the singular point is fixed by the involution $\sigma$, the curve $C_1$ is singular and $C_2$ is smooth. In particular, they are non-isomorphic and $C_1$, resp. $\tC_1$ is non-isomorphic to $X$, resp. $\tX$. Hence, the three curves and their double covers are also non-isomorphic for a generic choice of $\rho_4$.
\end{proof}

\section{\bf The Clifford Indices of the counterexamples}

In order to give an estimate for the Clifford index of the curve $X$ we use the following well known 
inequality of Castelnuovo (see \cite{c}):
Suppose $X$ is a smooth projective curve admitting two maps $f_1: X \ra Y_1$, $f_2 : X\ra
Y_2$ of degrees $n_1, n_2 \geq 2$ which do not both factor via a 
map $X \ra Z$ of degree $\geq 2$. Then
\begin{equation} \label{eq3.1}
g_X \leq (n_1-1)(n_2-1) + n_1g_{Y_1} + n_2g_{Y_2}.
\end{equation} 
For a modern proof of this inequality, note that the assumption implies that the map $(f_1,f_2): X \ra Y_1 \times Y_2$ is birational onto its image 
and apply the adjunction formula.

Moreover, if $\cliff X$ and $\gon X$ denote the Clifford index and the gonality of $X$, recall the following
inequality, which is valid for any smooth projective curve (see \cite{elms}).
\begin{equation} \label{eq3.2}
\cliff X + 2 \leq \gon X \leq \cliff X + 3.
\end{equation}
Finally recall that a covering is called {\it simple}, if it cannot be written as a composition of 2 coverings of degree $\geq 2$. 
\begin{lemma} \label{lem3.3}
Let $f: X \ra Y$ be a simple covering of degree $n$ of smooth projective curves with ramification divisor of degree $\delta$.
If $\delta \geq 2(n-1)n \gon Y$, then
$$
\gon X = n \gon Y.
$$

\end{lemma}
\begin{proof}
Certainly we have $\gon X \leq n \gon Y$. Suppose that $g:X \ra \bP^1$ is a map of degree $m < n \gon Y$. 
Since $f$ is a simple covering, we may apply inequlity \eqref{eq3.1} to $f$ and $g$, which gives, using the Hurwitz formula,
$$
n(g_Y -1) + 1 + \frac{\delta}{2} = g_X \leq (n-1)(m-1) + ng_Y.
$$
This implies
$$
\delta \leq 2(n-1)m < 2 (n-1) n \gon Y,
$$
which contradicts the assumption.
\end{proof}

\begin{corollary} \label{cor3.4}
If $n = 4$ and $\rho_4: X \ra Y$ is a simple covering of a general curve $Y$ of genus
$g_Y$ with ramification divisor of degree $\delta \geq 24 \gon Y$.
Then
$$
\cliff X \geq 2g_Y-1.
$$
\end{corollary}

\begin{proof}
For a general curve $Y$ of genus $g_Y$ we have $\gon Y = [\frac{g+3}{2}]$. Now Lemma \ref{lem3.3} and inequality \eqref{eq3.2} give
$$
\cliff X \geq \gon X -3 = 4 \gon Y -3 \geq 2g_Y - 1.
$$
\end{proof}

Certainly this estimate is not the best possible. Moreover, one can use the same method to give a more precise result
for any curve $Y$ with given gonality or Clifford index.

\section{\bf The class of $\tC_i$ in the Prym variety $P$}

\subsection{}
Let the notation again be as in Section 1. In particular $\rho_n: X \ra Y$ is a simply ramified cover of degree $n \geq 3$. 

In this section we compute the class of the image of $\tC_i$ in the Prym variety $P = \Prym (\kappa)$.

\subsection{}
We map $\tC_i$ into the Prym variety $P$ of $\kappa: \tX \lra X$ via the composition 
\[
\tC_i \lra J\tX \lra J\tX/\kappa^*JX = P.
\]
We shall compute the cohomology class $[\tC_i]$ of the image of $\tC_i$ by this map via degeneration to the case 
where $Y$ is an irreducible rational curve of arithmetic genus $g_Y$.  First we note

\begin{proposition}\label{propsameclass}
Assume $n\geq 3$. The curves $\tC_1$ and $\tC_2$ have the same cohomology class in
$P$.
\end{proposition}

\begin{proof}
The proof of \cite[Proposition 1 page 360]{beauville82} goes through
without change.
\end{proof}

In the case where $Y \simeq \bP^1$ the cohomology class  $[\tC_i]$ in $P$ was computed by Beauville 
in the proof of \cite[Proposition 2, p. 363]{beauville82} to be
\[
[\tC_i] = 2^{n -1} \frac{[\T_P]^{(g_X - 2)}}{(g_X-2)!}
\]
where $\T_P$ is the principal polarization of $P$ as the Prym variety of $\kappa$.
Here we generalize this formula to the case where $Y$ is not isomorphic to $\bP^1$.

\begin{remark}

In \cite[Proposition 2, p. 363]{beauville82}, Beauville could further divide his class by $4$ because his map from $\tC$ into $P$ factored through $2$-multiplication $2_P : P\ra P$. Except when $n=4$, this is not necessarily the case in our situation and we cannot divide the class by $4$.

\end{remark}

We first generalize the formula to the case where $Y \cong \bP^1$ but the double
cover $\tX \ra X$ is ramified. Note that in the ramified case, the curve $\tC$ does not
necessarily split into the union of the two curves $\tC_1$ and $\tC_2$. It is however defined in the same way as a subscheme of $\tX^{(n)}$ by the Cartesian diagram

\begin{eqnarray} \label{eq8}
\tC & \lra & \tX^{ (n) } \nonumber \\
\downarrow & & \downarrow \kappa^{ (n) } \\
Y & \lra & X^{ (n) }. \nonumber 
\end{eqnarray}

It has therefore a well-defined cohomology class in $\tX^{(n)}$. We push this class forward to the jacobian $J\tX$ and then to the Prym variety $P = J\tX / JX$. We denote this push-forward class by $[\tC]$: it is a well-defined element of the cohomology ring of $P$.

\begin{proposition}
Consider coverings
\[
\tX \stackrel{\kappa}{\lra} X \stackrel{\rho_n}{\lra} Y \cong \bP^1
\]
where $\kappa$ is a ramified double cover and $\rho_n$ is a ramified $n$-sheeted cover. Define the 
curve $\tC$ as before and map $\tC$ to the Prym variety $P$ of $\kappa$.
Then the push-forward class $[\tC]$ is
\[
[\tC] = 2^{n - g_X+1} \frac{[\T_P]^{(g_{\tX} - g_X - 1)}}{(g_{\tX} - g_X-1)!}.
\]
\end{proposition}

\begin{proof}
Let $\varphi_{X,n} : X^{(n)} \ra JX$ and  $\varphi_{\tX,n} : \tX^{(n)} \ra J\tX$
denote Abel maps such that the following diagram commutes
\[
\begin{array}{ccc}
\tX^{(n)} & \stackrel{\varphi_{\tX,n}}{\lra} & J\tX \\
\kappa^{(n)} \downarrow & & N \downarrow\\
X^{(n)} & \stackrel{\varphi_{X,n}}{\lra} & JX
\end{array}
\] 
where $N$ denotes the norm map. Let $\eta_{X,n}$ and $\eta_{\tX,n}$ be the respective 
cohomology classes of $X^{(n-1)} + p \subset X^{(n)}$ and 
$\tX^{(n-1)} + \widetilde{p} \subset \tX^{(n)}$ where $p$ and $\widetilde{p}$ are points of $X$ and $\tX$. Hence 
\begin{equation} \label{eqnvarphi*eta}
{\varphi_{\tX,n}}_* \eta_{\tX,n}^k = \frac{[\T_{\tX}]^{k + g_{\tX} - n}}{(k + g_{\tX} -n)!}
\end{equation}
where $\T_{\tX}$ is the principal polarization of $\tX$ as the jacobian of $\tX$.
As in \cite[p. 363]{beauville82} we rewrite Macdonald's formula \cite[p. 337]{macdonald62} for the class of the image of $Y = g^1_n$ in $X^{(n)}$ as 
\[ 
[Y] = \sum_{\alpha = 0}^{n-1} {n-1-g_X  \choose \alpha} \eta_{X,n}^{\alpha} \frac{\varphi_{X,n}^*[\T_X]^{n-1- \alpha}}{(n-1- \alpha)!}.
\]
Therefore, in $\tX^{(n)}$, the cohomolgy class of $\tC$ is
\[
[\tC]_{\tX^{(n)}} = \sum_{\alpha = 0}^{n-1} {n-1-g_X  \choose \alpha}2^{\alpha} \eta_{\tX,n}^{\alpha} \frac{\varphi_{\tX,n}^*N^*[\T_X]^{n-1- \alpha}}{(n-1- \alpha)!}.
\]
Using \eqref{eqnvarphi*eta} and the projection formula the push-forward of this to $J\tX$ is
\[
[\tC]_{J\tX} = \sum_{\alpha = 0}^{n-1} {n-1-g_X  \choose \alpha}2^{\alpha} \frac{[\T_{\tX}]^{\alpha + g_{\tX} -n}}{(\alpha + g_{\tX} -n)!} \frac{N^*[\T_X]^{n-1- \alpha}}{(n-1- \alpha)!}.
\]
It is easy to see that it follows from \cite[pp. 329-330]{mumford74} that
\[
4 [\T_{\tX}] = 2 N^*[\T_X] + q_n^*[\T_P]
\]
where $q_n: J\tX \ra P$ is the projection. Inserting, developping and using the fact that
\[
\left\{ {q_n}_* N^* : H^k(JX, \bQ) \lra H^{k-2g_X}(P,\bQ) \right\}  = \left\{ \begin{array}{ll}
0 & \mbox{if}\; k \neq 2g_X \\
2^{2(g_{\tX} - g_X)} & \mbox{if}\; k = 2g_X,
\end{array} \right.
\] 
we obtain the class $[\tC]$ in $P$:
\[
[\tC] = 2^{n -g_X+1} \sum_{\alpha = 0}^{n-1} {n-1-g_X \choose \alpha} {g_X \choose n-1-\alpha} \frac{[\T_P]^{g_{\tX} -g_X -1}}{(g_{\tX} -g_X -1)!} = 2^{n -g_X+1} \frac{[\T_P]^{g_{\tX} -g_X -1}}{(g_{\tX} -g_X -1)!}
\] 
because $\sum_{\alpha = 0}^{n-1} {n-1-g_X \choose \alpha} {g_X \choose n-1-\alpha} = 1$ (see \cite[p. 364]{beauville82}).
 \end{proof}

\subsection{} \label{subsecsing}
Now we generalize the construction to the case where $Y$ is an irreducible rational nodal curve 
of arithmetic genus $g_Y$, $\rho_n: X \ra Y$ is a simply ramified $n$-sheeted cover with branch 
locus disjoint from the nodes of $Y$, $\kappa: \tX \ra X$ is a $2$-sheeted cover, ramified at 
the nodes of $X$ such that near each ramification point the covering involution does not 
exchange the two branches of $\tX$. Note that the latter is a type  of Beauville admissible 
double cover. Consider the Cartesian diagram of normalizations:

\begin{equation*}
\begin{xy}
\xymatrix@R=24pt@C=24pt{
\tX_{\nu} \ar[d] \ar[r] & \tX \ar[d]\\
X_{\nu} \ar[d] \ar[r] & X \ar[d]\\
\bP^1 \ar[r] & Y}
\end{xy}
\end{equation*}
which induces the Cartesian diagram (the curves $\tC$ and $\tD$ below are defined by the diagram)

\begin{equation} \label{diag1.1}
\begin{xy}
\xymatrix@R=16pt@C=5pt{   
      &  \tX_{\nu}^{(n)}  \ar[rrrr] \ar[dd] &&&& \tX^{(n)} \ar[dd]\\
      \tD \ar@{^{(}->}[ur] \ar[rrrr] \ar[dd] &&&& \tC \ar@{^{(}->}[ur] \ar[dd] \\
      &  X_{\nu}^{(n)} \ar[rrrr]   & &&& X^{(n)}\\
      \bP^1 \ar[ur] \ar[rrrr] &&&& Y \ar[ur] &   }\end{xy}
\end{equation}
and hence the diagram

\begin{equation*}
\begin{xy}
\xymatrix@R=24pt@C=24pt{
\tC \ar[d] & \tD \ar[d] \ar[l]\\
J\tX \ar[d] \ar[r] & J\tX_{\nu} \ar[d]\\
P \ar[r]^{\alpha} & P_{\nu}
}\end{xy}
\end{equation*}
where $P$ and $P_{\nu}$ are, respectively, the generalized Prym varieties of the covers $\tX \ra X$ and $\tX_{\nu} \ra X_{\nu}$. 

Let $\beta: P_{\nu} \ra P$ be the pseudo-inverse of $\alpha$, i.e., $\beta$ is the unique isogeny such that $\alpha \circ \beta = 2_{P_{\nu}}$ and $\beta \circ \alpha = 2_P$. Note that the degree of $\alpha$ is $2^{n g_Y -1}$. Hence the degree of $\beta$ is $2^{2g_X -ng_Y -1}$.

\begin{proposition}
We have
\[
\beta_*\frac{[\T_{P_{\nu}}]^{g_X-2}}{(g_X - 2)!} = 2^{g_X - ng_Y + 1}\frac{[\T_{P}]^{g_X -2}}{(g_X -2)!}.\
\]
where $\T_{P_{\nu}}$ is the principal polarization of $P_{\nu}$ as the Prym variety of the double cover $\tX_{\nu}\ra X_{\nu}$.
\end{proposition}

\begin{proof}
By \cite[p. 159]{beauville771} the map $\alpha$ is an isogeny with kernel an isotropic subgroup with respect to $\T_{\tX}|_P$ of points of order 2. This implies that $\alpha^*[\T_{P_{\nu}}] = [\T_{\tX}|_P]$, equivalently, $\alpha^*[\T_{P_{\nu}}] = 2[\T_P]$. 
The latter implies $\beta^*[\T_{P}] = 2[\T_{P_{\nu}}]$. Therefore 
\[
\beta^*\frac{[\T_{P}]^{g_X -2}}{(g_X -2)!} = 2^{g_X-2}\frac{[\T_{P_{\nu}}]^{g_X-2}}{(g_X - 2)!}.
\]
Poincar\'e duality gives the following commutative diagram
\[
\begin{array}{ccc}
H^{2g_X -4}(P, \bQ) & \stackrel{\beta^*}{\lra} & H^{2g_X -4}(P_{\nu}, \bQ)\\
d \downarrow && \downarrow\\
H_2(P, \bQ) & \stackrel{\beta_*}{\lla} &H_2(P_{\nu}, \bQ)
\end{array}
\] 
where the right hand vertical map is Poincar\'e duality and the left hand vertical map is Poincar\'e duality multiplied by the degree $d = 2^{2g_X -n g_Y -1}$ of $\beta$. Now the proposition follows.
\end{proof} 

\begin{corollary}
The class $[\tC]$ in $P$ is 
\[
[\tC] = 2^{n}\frac{[\T_{P}]^{g_X -2}}{(g_X -2)!}.
\]
\end{corollary}

\begin{proof}
Using the previous two propositions and the fact that $\beta \circ \alpha = 2_P$, we compute
\[
[\tC] = \frac{1}{4}\beta_* [\tD] = \frac{1}{4} 2^{n-g_{X_{\nu}}+1} \beta_* \frac{[\T_{P_{\nu}}]^{g_X-2}}{(g_X-2)!} = 
2^{n - g_{X_{\nu}} -1 + g_X -ng_Y +1} \frac{[\T_P]^{g_X -2}}{(g_X - 2)!}.
\]
\end{proof}

\begin{theorem} \label{thmclass}
For any unramified double cover $\kappa: \tX \ra X$ and any simply ramified $n$-sheeted cover $\rho_n: X \ra Y$ as in 
Section 1 the class of the image of $\tC_i$ in the Prym variety $P$ of $\kappa$ is 
\[
[\tC_i] = 2^{n - 1}\frac{[\T_{P}]^{g_X -2}}{(g_X -2)!}.
\] 
\end{theorem}

\begin{proof}
The class $[\tC]$ of the image of $\tC$ can be computed from the corollary by a simple degeneration argument using the fact that the constructions can be 
done in continuous families including admissible covers where the Prym varieties do not degenerate and the 
fact that the above cohomology classes live in a local system of free abelian groups. We
then use Proposition \ref{propsameclass} to obtain the classes $[\tC_i]$.
\end{proof}

\bibliographystyle{plain}


\end{document}